\documentclass[conference]{IEEEtran}
\IEEEoverridecommandlockouts
\newcommand{\R}{\mathbb{R}}
\newcommand{\ep}{\varepsilon}
\newcommand{\F}{\mathcal{F}}

\newcommand{\B}{\mathcal{B}}
\newcommand{\PO}{\mathcal{P}}
\newcommand{\X}{\mathcal{X}}
\newcommand{\m}{\tilde{\mu}}
\newcommand{\W}{\overline{W}}

\usepackage{cite}
\usepackage{amsmath,amssymb,amsfonts}
\usepackage{algorithmic}
\usepackage{graphicx}
\usepackage{textcomp}
\usepackage{xcolor}
\usepackage{mathrsfs}
\usepackage{latexsym}
\newtheorem{theorem}{Theorem}[section]
\newtheorem{lemma}[theorem]{Lemma}

\newtheorem{corollary}[theorem]{Corollary}

\numberwithin{equation}{section}
\def\BibTeX{{\rm B\kern-.05em{\sc i\kern-.025em b}\kern-.08em
    T\kern-.1667em\lower.7ex\hbox{E}\kern-.125emX}}
\begin{document}
\title{A theoretical basis for model collapse in recursive training
\thanks{Research supported by a grant from Google Research Asia}}

\author{\IEEEauthorblockN{Vivek S.\ Borkar}
\IEEEauthorblockA{\textit{Department of Electrical Engineering (retd.)} \\
\textit{Indian Institute of Technology Bombay}\\
Mumbai 400076, India \\
Email: borkar.vs@gmail.com}
}

\maketitle

\begin{abstract}
It is known that recursive training from generative models can lead to the so called `collapse' of the simulated probability distribution. This note shows that one in fact gets two different asymptotic behaviours depending on whether an external source, howsoever minor, is also contributing samples.
\end{abstract}


\section{Introduction} Several recent works have observed and discussed a phenomenon dubbed `model collapse' in recursive training by generated data where, at each stage, the present data is used as the seed for generating new data for the next stage. See, e.g., \cite{Alem, Elvis1, Elvis2, Gerst, Marchi, Suresh, Shum1, Shum2}. The objective of this note is to give a precise mathematical formulation of the problem and explain the phenomenon.  This work depends heavily on certain tools of advanced probability, notably the so called `weak convergence' of probability measures, discrete parameter martingales, and Markov chains on a Polish space. The minimal necessary background for our purposes is recalled in the next section.
Sections 3 and 4 establish our main results for the simplest of the mechanisms for recursive training. The first case, considered in Section 3, has what we call `\textit{persistent excitation}', borrowing a phrase from control theory. Here, there is  injection of fresh data at every step over and above the recursive generation of data. This leads to a phenomenon which we call `\textit{degeneration}' in order to distinguish it from collapse. The reason for doing so will be self-evident. Section 4 considers the case without continued injection of fresh data. This leads to a collapse in a sense that we make precise and is in the spirit of the common notions of collapse in literature. An appendix recalls a key result from \cite{BorkarTopics} used in the proof.\\

We work with samples in a general Polish space, defined in the next section, as opposed to taking the simpler option of either a finite set or a euclidean space. This is in order to accommodate function spaces such as $C([0,T]; \R^d) :=$ the Banach space of continuous functions $[0,T] \mapsto \R^d$ with norm $\|f\| := \max_{t\in[0,T]}\|f(t)\|$. This does arise in the context of diffusion models and finance. More importantly, Polish spaces also have the property that the space of probability measures on a Polish space with the `topology of weak convergence' is itself a Polish space, a fact we use crucially. (This is not true in general, e.g., for Banach spaces - the space of probability measures on a Banach space with the above topology is \textit{not} a Banach space.) Also, Polish spaces are at the right level of generality where bulk of the results of standard probability theory apply without having to weaken them.

\section{Mathematical background}

Our analysis will draw heavily upon the three topics in probability theory mentioned above.  We briefly summarize the relevant results here. These can be found respectively, in \cite{Bill} (see also \cite{KRP} for a more extensive treatment), \cite{Neveu}, and \cite{Benaim} (see also \cite{Meyn} for a more extensive treatment),  respectively.\\

\subsection{Convergence of probability measures :}\label{convmeas} 

Let $S$ be a Polish space, i.e., a separable topological space with its topology compatible with a complete metric. Let $\B$ denote its Borel $\sigma$-field, i.e., the smallest $\sigma$-field containing its open sets. Let $C_b(S) :=$ the space of bounded continuous functions $S \to \mathbb{R}$. Let $\PO(S)$ denote the space of probability measures on $(S, \B)$, with the coarsest topology required to render continuous the maps $\mu \in \PO(S) \mapsto \mu(f) := \int fd\mu \in \mathbb{R}$ for any $f \in C_b(S)$. Then $\PO(S)$ is also a Polish space with one choice of a complete metric being given by the Prohorov metric. See \cite{Bill}, Appendix III, for details. Likewise, $\PO(\PO(S))$ is a Polish space, a fact we use later. 
%If $\mu_n(f) \to \mu(f) \ \forall f \in C_b(S)$ with compact supports, but $\mu$ is not a probability measure, this corresponds to the notion of \textit{vague} convergence (equivalently, convergence in the \textit{vague} topology) defined for finite non-negative measures on $S$. 

We shall use the following definitions and results from \cite{Bill}:

\begin{enumerate}

\item A set $A \subset \PO(S)$ is said to be \textit{tight} if given any $\ep > 0$, there exists a compact set $K_\ep \subset S$ such that $\mu(K_\ep) > 1 - \ep$ for all $\mu \in A$. Then a set $B \subset \PO(S)$ is relatively compact (equivalently, relatively sequentially compact) if and only if it is tight (Prohorov's theorem, \cite{Bill}, p.\ 37).  A singleton, i.e., $\{\mu\}, \mu \in \PO(S)$, is always tight by a theorem of Oxtoby and Ulam (\cite{Bill}, p.\ 10). See \cite{Bill}, Appendix III, for details.\\

\item There is a countable collection $\{f_i, i \geq 1\} \subset C_b(S)$ such that
$$\mu_n \to \mu \ \mbox{in} \ \PO(S) \ \Longleftrightarrow \ \mu_n(f_i) \to \mu(f_i) \ \forall i.$$
Such a collection $\{f_i\}$ is known as a \textit{convergence determining class} \cite{Bill}. By adding suitable positive constants to the $f_i$'s  and then rescaling them by suitable positive scalars if necessary, we may assume without loss of generality that $0 \leq f_i(\cdot) \leq 1 \ \forall \  i \geq 1$. It is easy to see that the $\{f_i\}$ separate points of $\PO(S)$, i.e., $\mu(f_i) = \nu(f_i) \ \forall i$ implies $\mu = \nu$. See \cite{Bill}, p.\ 15. \\

\item Let $A \subset \PO(S)$. For $\mu \in \PO(S)$, let $\zeta_\mu$ be a $\PO(S)$-valued random variable such that for all Borel $B \subset S$,
$$\mu(B) = E[\zeta_\mu(B)] \ \mbox{for all} \ \mu \in A. \ \ \ \ \ \ \ \ \  (\dagger)$$
 Let $\X(\mu) \in \PO(\PO(S))$ denote the law of $\zeta_\mu$ for $\mu \in A$. Then $A$ is tight in $\PO(S)$ if and only if the set 
$\bar{A} := \{\X(\mu), \mu \in A\},$ is tight as a subset of  $\PO(\PO(S))$. This is Lemma 2.1, p.\ 106, of \cite{BorkarTopics}. This result and its proof are recalled in the Appendix as Theorem \ref{barycenter}.  
Whenever $(\dagger)$ holds, we shall say that $\mu$ is a \textit{barycenter} of $\X(\mu)$ or, by abuse of notation, of $\zeta_\mu$. Likewise, a $\PO(S)$-valued random variable $\mu'$ is said to be a \textit{conditional} barycenter of $\zeta$ or the law thereof, if  $\mu'(B) = E[\zeta(B)|\F']$ a.s.\ for all Borel $B \in \F$ and some sub-$\sigma$-field $\F'$ of $\F$.\\
\end{enumerate}

\subsection{Martingales :}\label{martin} 
Consider a probability space $(\Omega,\F,P)$ endowed with an increasing sequence of sub-$\sigma$-fields $\F_n, n \geq 0$, of $\F$ (i.e., $\F_n \subset \F_{n+1} \subset \F \ \forall  \ n \geq 0$). A sequence of $\R^d$-valued random variables $M(n), n \geq 0$, satisfying  
$$E[\|M(n+1)\| | \F_n] < \infty \ \forall \ n,$$ 
is said to be a martingale with respect to $\{\F_n\}$ (equivalently, an  $\{\F_n\}$-martingale) if $\forall n$, $M(n)$ is $\F_n$-measurable and $E[M(n+1)|\F_n] = M(n).$ 
If $\sup_n\|M(n)\| \leq K^* < \infty$ a.s.\ for a constant $K^*$ (i.e., the martingale is uniformly bounded by a constant, a.s.), then $M(n) \to$ a random variable  $M(\infty)$ in $\R^d$, both a.s.\  and in mean. Furthermore, $M(n) = E[M(\infty)|\F_n]$ a.s.\ $\forall n$. This is a special case of a much more general result (See Proposition IV-2-3, p.\ 65, \cite{Neveu}.) The version above will suffice for our purposes.

\subsection{Markov chains :} 
An $S$-valued process $X_n, n \geq 0,$ is a (time-homogeneous) Markov chain on state space $S$ with a transition probability kernel 
$$x \in S \ \mapsto \ p(dy|x) \in \PO(S)$$ 
if, for $n \geq 0$,
$$P(X_{n+1} \in A | X_m, m \leq n) = p(A|X_n) \ \  \forall \ A \in \B,$$
or equivalently, if
$$E[f(X_{n+1}) | X_m, m \leq n] = \int f(y)p(dy|X_n) \ \ \forall \ f \in C_b(S).$$
We assume that the map  $x \in S \mapsto p(dy|x) \in \PO(S)$ is continuous. (This makes the Markov chain a \textit{Feller} process.) This Markov chain is said to be \textit{stable} if for any initial law $\mu_0$ of $X_0$, the laws of $X_n, n \geq 0,$ remain tight in $\PO(S)$. If so, then the Cesaro sums 
$$\frac{1}{n}\sum_{m=0}^{n-1}E[f(X_m)] \to \int fd\nu  \ \ \forall  \ f \in C_b(S)$$
for a $\nu$ that is an invariant probability measure for the Markov chain. That is,
$$\nu(A) = \int p(A|x)\nu(dx) \ \ \forall \ \mbox{Borel} \  A \in \B$$
$$\left(\Longleftrightarrow \ \int fd\nu = \int\int f(y)p(dy|x)\nu(dx) \ \ \forall \ f \in C_b(S)\right).$$
This $\nu$ need not be unique and will in general depend on the initial condition $X_0$. \\

Under our assumptions of stability and continuity of $p(dy| \cdot )$, the set of invariant measures is nonempty, closed and convex. Furthermore, its extreme points are mutually singular, i.e., are concentrated on a collection of sets that form a partition of  a subset of $S$, or of $S$ itself if there are no transient states. This is known as the Doeblin decomposition. See \cite{Meyn} for this and much more.

\section{Degeneration in presence of persistent excitation}

We  analyze a basic model for recursive training based on the following scenario. Fix $N \gg 1$. We begin with a probability measure $\mu_0 \in \PO(S)$ and  a parametrized family $P_\theta$ for  $\theta \in$ some convex compact parameter set $D \subset \R^p$ (say) containing the `true' parameter  $\theta^*$ such that $P_{\theta^*} = \mu_0$. We assume that the map $\theta \in D \mapsto P_\theta \in \PO(S)$ is continuous. Then, since $D$ is compact, so will be the set $\{P_\theta, \theta \in D\}$. In particular, it will be tight. \\

We treat $\mu_0$ as the distribution of i.i.d.\ samples generated by an external source such as an experimental set-up or physical measurement.  Let $a \in [0,1]$ and $b, c \in [0,1]$ with $b+c=1$. The aforementioned  samples are  augmented by samples with distribution $\nu_n := bP_{\theta(n)} + c\mu_n$ generated by a generative model at time $n$, see details below.  Let $\theta(0)$ have  some prescribed distribution with barycenter $\mu_0$ . Recursively  do the following for $n \geq 0$:\\

\begin{enumerate}

\item Let $\nu_n := bP_{\theta(n)} + c\mu_n.$\\

\item Generate $N$ independent samples $X^n_1, \cdots , X^n_N$ according to $\mu_0$ with probability $a$ and according to $\nu_n$ with probability $1-a$ (i.e., according to $P_{\theta(n)}$ with probability $(1-a)b$ and according to $\mu_n$ with probability $(1-a)c$).\\

\item Recursively generate $\theta_{n+1}$ from $\theta_n$ and $\{X^n_i\}$, and generate the corresponding $P_{\theta_{n+1}}$ such that $P_{\theta_n}$ is the conditional barycenter of $P_{\theta_{n+1}}$ in the sense that $$E[P_{\theta_{n+1}}(A)|\F_n] = E[P_{\theta_{n+1}}(A)|\theta_n] = P_{\theta_n}(A) \ \mbox{a.s.}$$
for Borel $A \subset S$. This will follow, e.g., if the law of $\theta(n)$ is the conditional barycenter of the law of $\theta(n+1)$.   (In practice, this may be an idealization, not an `exact fit'.) We also assume that the conditional law of $\theta(n+1)$ given $\theta(n)$ is continuous as a map $D \mapsto \PO(D)$ and is independent of $n$.\\

\item Compute the empirical distribution $\mu_{n+1}$ of the samples according to
$$\mu_{n+1}(A) := \frac{1}{N}\sum_{i=1}^NI\{X^n_i \in A\}, \ A \in \B,$$
where the `indicator function' $I\{ \cdots \} = 1$ if `$\cdots$' is true and $0$ otherwise.\\

\item $n \to n+1$.
\end{enumerate}

\medskip

 For the purpose of our analysis, we may treat $bP_{\theta(n)} + c\mu_n$ as a single measure that is not necessarily finitely supported, denoted henceforth as $\breve{\mu}_n$ by abuse of terminology. We have the following.\\

\begin{lemma}\label{tightness} The laws of $\breve{\mu}_n, n \geq 0,$ are tight in $\PO(\PO(S))$.
\end{lemma}

\medskip

\noindent \textit{Proof} Letting $\check{\mu}_n \in \PO(S)$ denote the barycenter of $\breve{\mu}_n$, it follows from our construction that\
$$\check{\mu}_{n+1}(f) = a\mu_0(f) + (1-a)\check{\mu}_n(f) \ \ \forall \ f \in C_b(S).$$
Thus $\check{\mu}_n(f) \to \mu_0(f)$, implying in particular that $\{\check{\mu}_n\}$ is a tight set. The claim now follows from Theorem \ref{barycenter} in the Appendix. \hfill $\Box$

\medskip

\begin{lemma} $\{(\mu_n, \theta_n)\}$ is a $\PO(S)\times D$-valued time-homogeneous, stable Markov chain with a continuous transition kernel. \end{lemma}

\medskip

\noindent \textit{Proof} By construction, the conditional law  of $(\mu_{n+1}, \theta_{n+1})$ given $(\mu_m, \theta_m), m \leq n,$ is the same as its conditional law given $(\mu_n, \theta_n)$ alone, for $n \geq 0$. Furthermore, it is clearly independent of $n$. This proves that it is a time-homogeneous Markov chain. Denote the aforementioned  conditional law as $\kappa( \cdot |(\mu_n,\theta_n))$, which is then the transition kernel for the Markov chain $\{(\mu_n,\theta_n)\}$. Stability follows from Lemma \ref{tightness}. To prove the continuity of the transition kernel, we need to verify that
$\int \varphi(\mu',\theta')\kappa(d(\mu', \theta')|\mu,\theta)$ is continuous in $(\mu, \theta)$ for all $\varphi \in C_b(\PO(S)\times. D)$. Let $\hat{\nu}_n \to \hat{\nu}_\infty$ in $\PO(S)\times D$. Then $ L := \{\hat{\nu}_n, 1 \leq n \leq \infty\}$ is compact in $\PO(S)\times D$. Consider maps $\varphi \in C_b(\PO(S)\times D)$ of the form
\begin{eqnarray*}
\lefteqn{(\mu, \theta) \in \PO(S)\times D \mapsto
\varphi((\mu, \theta)) } \\
&=& \Psi\left(\int g_1d\mu, \cdots , \int g_md\mu, \theta\right) \in \R
\end{eqnarray*}
for some $m \geq 1$, $g_i \in C_b(S), 1 \leq i \leq m, \Psi \in C_b(\R^m\times D)$. By the Stone-Weierstrass theorem, such functions restricted to $L$ are dense in $C(L)$. Hence it suffices to verify that for $\varphi$ as above,
$$\int \varphi(\nu')\kappa(d\nu'|\hat{\nu}_n) \to \int \varphi(\nu')\kappa(d\nu'|\hat{\nu}_\infty).$$
This follows a.s.\ by direct verification from the construction of $\nu'$ given $\hat{\nu}_{(\cdot)}$, which is identical to the construction of $(\mu_{n+1},\theta_{n+1})$ from $(\mu_n, \theta_n)$.  
%Specifically, we have, for $\{X^i_k, 1 \leq i \leq m, 1 \leq k \leq N\}$ i.i.d.\ with law $\kappa( \cdot | \nu )$,
%\begin{eqnarray*}
%\lefteqn{\int \varphi(\nu')\kappa(d\nu'|\hat{\nu}_n)} \\
%&=& E\left[\Psi\left(\frac{1}{N}\sum_{k=1}^Ng_1(X^1_k), \cdots ,  \frac{1}{N}\sum_{k=1}^Ng_m(X^m_k)\right) \Big| \hat{\nu}_n\right] \\
%&=& \int \Psi\left(\frac{1}{N}\sum_{k=1}^Ng_1(x^1_k), \cdots ,  \frac{1}{N}\sum_{k=1}^Ng_m(x^m_k)\right)(d\hat{\nu}_n)^{mN} \\
%&\to& \int \Psi\left(\frac{1}{N}\sum_{k=1}^Ng_1(x^1_k), \cdots ,  \frac{1}{N}\sum_{k=1}^Ng_m(x^m_k)\right)(d\hat{\nu}_\infty)^{mN} \\
%&\to& \int \varphi(\nu')d(\kappa(\nu'|\hat{\nu}_\infty),
%\end{eqnarray*} 
%which proves the claim. 
\hfill $\Box$

\medskip

This implies in particular that the set of stationary distributions for this Markov chain is nonempty closed and convex with its extreme points corresponding to the ergodic measures. This leads to our first main result.\\

\begin{theorem} If $a > 0$, the barycenter of the stationary distributions of the $\PO(S)\times D$-valued Markov chain $\{(\mu_n, \theta_n)\}$ is $(\mu_0, \theta_0)$. \end{theorem}

\medskip

\noindent \textit{Proof} It is easy to see that the barycenter $(\bar{\mu}, \bar{\theta})$ of any stationary distribution of the Markov chain $\{(\mu_n, \theta_n)\}$ above must be a fixed point of the map $\nu \mapsto a(\mu_0, \theta_0) + (1-a)\nu$.  This implies that $\bar{\mu} = \mu_0, \bar{\theta} = \theta_0$. \hfill $\Box$

\medskip

Note that this does not imply the uniqueness of the stationary distribution. Consider, e.g., uniform distributions on the sets of Dirac measures on the points on concentric circles of various radii. Their barycenter is the Dirac measure at their common centre, though their supports are disjoint.

What we have, in comparison with i.i.d.\ samples generated from $\mu_0$ alone, is the samples generated according to random measures $\{\mu_n\}$ whose barycenter is $\mu_0$. There is degradation in quality on two counts. First, the samples are no longer i.i.d., $\{\mu_n\}$ being a Markov chain. More importantly,  while the mean behaviour is maintained because $\int fd\mu_0 = E[\int fd\mu_n]$ for all $f \in C_b(S)$, the new samples are more `noisy' in the sense that
$E[g(\int fd\mu_n)] \geq g(\int fd\mu_0)$ for all convex $g: \R \mapsto \R$ by conditional Jensen's inequality. The extent thereof will depend on the magnitude of the parameter $a$. We call this phenomenon as `\textit{degeneration}' of sample quality in order to differentiate it from \textit{collapse} which we study next.\\

\section{Total collapse in absence of persistent excitation}

Next we consider the case when $a = 0$. Let 
$$\F_n := \sigma(X^m_i, 1 \leq i \leq N, 0 \leq m \leq n), \ n \geq 0.$$ 
Let $\{f_i\}$ be a countable convergence determining class for $\PO(S)$ with $0 \leq f_i(\cdot) \leq 1$. Let $M_i(n) = \mu_n(f_i), n \geq 0, i \geq 1$.\\

\begin{lemma}\label{martconv} For each $i \geq 1$, $M_i(n), n \geq 0,$ is an $\{\F_n\}$-martingale which converges a.s.\ and in mean to a random variable $M_i(\infty)$. Furthermore, $M_i(n) = E[M_i(\infty)|\F_n] \ \forall \ i, n$.\end{lemma}

\medskip

\noindent \textit{Proof} The martingale property follows from the fact that
$$E\left[\int f_id\mu_{n+1} \Big| \F_n\right] = \int f_id\mu_n \ \ \forall \ i, n.$$
Since the $|f_i|$ are bounded by $1$, so are the $\{|M_i(n)|, n \geq 0\}$ for each $i$. Hence they converge a.s.\ and in mean to limiting $\R^d$-valued random variables $M_i(\infty)$ respectively.  The final claim is immediate from the results recalled in Subsection \ref{martin}. \hfill $\Box$

\medskip

In particular, $\{M_i(n)\}$ can be identified with $(\R^d)^\infty$-valued random variables $[M_1(n), M_2(n), \cdots ]$. These then converge jointly, a.s.\ and in component-wise mean, to the  limiting random variable $[M_1(\infty), M_2(\infty), \cdots]$. 
Together with Lemma \ref{tightness}, this has the following important consequence.\\

\begin{corollary}\label{convergence} $\mu_n \to \mu_\infty$ a.s. for some $\PO(S)$-valued random variable $\mu_\infty$.
\end{corollary}

\medskip

\noindent \textit{Proof} Since the laws of $\{\mu_n\}$ are tight, they are sequentially relatively compact in $\PO(\PO(S))$. Let $\mu_*$ be a limit point in law of $\mu_n$ as $n\to\infty$, i.e., the laws of $\{\mu_n\}$ converge to the law of $\mu_*$ along a subsequence. Then along this subsequence,
$$\left[\int f_1d\mu_n, \int f_2d\mu_n, \cdots , \int f_id\mu_n, \cdots \right] \to$$
$$\left[\int f_1d\mu_*, \int f_2d\mu_*, \cdots , \int f_id\mu_*, \cdots \right]$$
in law. But by Lemma \ref{martconv}, the left hand side converges  to 
$[M_1(\infty), M_2(\infty), \cdots]$ a.s., hence in law, for $\{M_i(\infty)\}$ as defined in the proof thereof. It follows that the laws of $[M_1(\infty), M_2(\infty), \cdots]$ and $[\int f_1d\mu_*, \int f_2d\mu_*, \cdots ]$ agree. Letting 
$$\phi : \upsilon \in \PO(S) \mapsto \left[\int f_1d\upsilon, \int f_2d\upsilon, \cdots\right] \in [0,1]^\infty,$$ 
$\phi : \PO(S) \to \phi(\PO(S))$ is continuous, one-one and onto, because $\{f_i\}$ separates points of $\PO(S)$.  By the foregoing, $[M_1(\infty), M_2(\infty), \cdots] \in \phi(S)$ a.s. Define $\mu_\infty := \phi^{-1}\left([M_1(\infty), M_2(\infty), \cdots]\right)$. Then $\mu_n \to \mu_\infty$ a.s.\ in $\PO(S)$. Note also that $\mu_*$ above agrees with $\mu_\infty$ in law and therefore is the unique limit point in law of $\{\mu_n\}$. In other words, $\mu_\infty$ is the unique limit both in law  and a.s.\ of $\mu_n$ as $n\to\infty$. 
\hfill $\Box$

\medskip

In view of the properties of such Markov chains already noted above, the chain has a nonempty closed convex set of invariant probability measures. At the same time, Corollary \ref{convergence} tells us that $\mu_n \to \mu_\infty$ a.s. The only way a time-homogeneous stable Markov chain can converge a.s.\ is if its limiting invariant distributions are Dirac. That is, the chain converges to the set of its absorbing states, which are Dirac measures on $\PO(S)$. This leads to our main result:\\

\begin{theorem} $\mu_n \to \delta_\gamma$ a.s., where $\gamma$ is a $S$-valued random variable. \end{theorem}

\medskip

This is the precise statement (and proof) of model collapse in recursive training for a purely generative mechanism when the current data is the seed for generating the data for the next step.  It also coincides with the mechanism of \textit{stochastic localization} that has been analyzed in a different context, see, e.g., \cite{Montanari}\footnote{The author thanks Prof.\ Bruce Hajek for this observation.}.
Note also that  if $\mu_{n+1}$ is simply the empirical law of some fixed number (say, $N'$) of samples generated according to $\mu_n$ for $n \geq 0$, then the above convergence will be in a finite random  time.

\section{Discussion}

We have presented two distinct scenarios with completely different conclusions, depending on whether the samples are generated by a purely generative model feeding on itself, or there is also another genuine source of i.i.d.\ fresh data. In the former case, the measure-valued Markov chain on  $\PO(S)$  degenerates into a random Dirac measure on $\PO(S)$, in the latter they asymptotically form a stationary Markov chain (not necessarily unique) whose stationary distribution has its barycenter at the law of samples generated by the natural source.

What about the errors implicit in the idealizations underlying our stylized model? They will certainly matter, both in terms of perturbing the idealized time asymptotics described above and the actual finite time behaviour \textit{en route} to this asymptotics. Even otherwise, the rate of convergence \cite{Suresh}, finite time bounds, etc., remain important issues because they inform us about how many generative samples can be used for augmenting available data without seriously affecting its quality. These are problems for future research.

\section{Appendix}

\medskip

The following fact was used above. We recall its proof from \cite{BorkarTopics}, p.\ 106.  For $\mu \in \PO(S)$, let $\zeta_\mu$ be a $\PO(S)$-valued random variable such $\mu(B) = E[\zeta_\mu(B)]$ for all $B \subset S$ Borel. Let $\X(\mu) \in \PO(\PO(S))$ denote the law of $\zeta_\mu$ for $\mu \in \PO(S)$.

\medskip

\begin{theorem}\label{barycenter}  $A \subset \PO(S)$ is tight  if and only if the set $\tilde{A} := \{\X(\mu), \mu \in A\} \subset \PO(\PO(S))$ is tight. \end{theorem}

\medskip

\noindent \textit{Proof} Let  $A \subset \PO(S)$ and  
$$\tilde{A} := \{\X(\mu) \in \PO(\PO(S)) : \mu \in A\}.$$ 
Suppose $A$ is tight, but $\tilde{A}$ is not. Recall the map $\phi : \PO(S) \to [0,1]^\infty$ defined in the proof  of Corollary \ref{convergence}. Then $\phi : \PO(S) \to \phi(\PO(S))$ is seen to be continuous, one-one and onto. Let $\W := \overline{\phi(S)} \subset [0,1]^\infty$. $\overline{W}$ is compact and therefore can be viewed as a compactification of $W$. We identify $W$ with $\PO(S)$ via the homeomorphism $\PO(S) \longleftrightarrow \phi(\PO(S))$. Let $\partial W := \W\backslash W$.  Since $\tilde{A}$ is not tight in $\PO(\PO(S))$, we can find a sequence $\m_n, n \geq 1,$ in $\tilde{A}$ such that $\m_n \to \m_\infty \in \PO(\W)$ with $\m_\infty(\partial W) > 0$. Define the barycenters $\mu_n \in \W$ of $\m_n$ by 
$$\int fd\mu_n = E\left[\int fd\m_n\right], \ f \in C_b(\W), 1 \leq n \leq \infty.$$
Then $\{\mu_n, 1 \leq n < \infty\} \subset A \subset \PO(S),$ but $\mu_\infty \notin \PO(S)$. This contradicts the tightness of $A$. Thus $\tilde{A}$ must be tight in $\PO(\PO(S))$.\\

Conversely, suppose $Z \subset \PO(\PO(S))$ is tight. Let $\ep > 0$. Then by the definition of tightness, we can find  a compact set $K_\ep \in \PO(S)$ such that
$$\upsilon\left(K_\ep\right) \geq \ 1 - \frac{\ep}{2} \ \ \ \ \ \forall \ \upsilon \in Z.$$
Let $C \subset S$ be a compact set satisfying $\kappa(C) \geq 1 - \frac{\ep}{2}$ for all $\kappa \in K_\ep$. This is possible by Prohorov's theorem because $K_\ep$ is tight. Then for $\bar{\upsilon}$ defined by $\upsilon = \zeta_{\bar{\upsilon}}$ in our earlier notation,
\begin{eqnarray*}
\bar{\upsilon}(C) &=& E[\upsilon(C)] \\
&=& E[\upsilon(C)I\{\upsilon\in K_\ep\}]  + E[\upsilon(C)I\{\upsilon \notin K_\ep\}] \\
&\geq& \left(1 - \frac{\ep}{2}\right)^2 \\
&\geq& 1 - \ep.
\end{eqnarray*}
Therefore the set $A := \{\bar{\kappa} : \kappa \in \bar{A}\}$ is tight. \hfill $\Box$

\medskip

\section*{Acknowledgment}
The author  thanks Prof.\ Andrew Thangaraj for bringing this problem to his notice and Prof.\ K.\ S.\ Mallikarjuna Rao, Dr.\ Raghavendra Tripathi  and Satush Parikh for carefully reading the manuscript.

\end{document}